\documentclass[11pt,a4paper]{article}
\usepackage[top=2.5cm,bottom=2.5cm,left=2.5cm,right=2.5cm]{geometry}
\usepackage[utf8]{inputenc}
\usepackage{xcolor}
\usepackage{graphicx}
\usepackage{amsmath,amssymb,latexsym}
\usepackage{epsf,enumerate}
\usepackage{tcolorbox}
\usepackage{url}

\usepackage{bm}

\newtheorem{problem}{Problem}[section]
\newtheorem{conjecture}{Conjecture}[section]
\newtheorem{theorem}{Theorem}[section]
\newtheorem{proposition}{Proposition}[section]

\newenvironment{boxedproblem}{\begin{problem}}{\end{problem}}
\newenvironment{boxedconjecture}{\begin{conjecture}}{\end{conjecture}}
\newenvironment{boxedtheorem}{\begin{theorem}}{\end{theorem}}
\newenvironment{boxedproposition}{\begin{proposition}}{\end{proposition}}

\tcolorboxenvironment{boxedproblem}{}
\tcolorboxenvironment{boxedconjecture}{}
\tcolorboxenvironment{boxedtheorem}{}
\tcolorboxenvironment{boxedproposition}{}

\usepackage[super]{nth}

\usepackage{hyperref}

\makeatletter
\let\Hy@linktoc\Hy@linktoc@none
\makeatother

\title{Open problems of the \nth{33} Workshop on Cycles and Colourings}

\author{(edited by Alfréd Onderko)}

\date{}

\begin{document}

\maketitle

\begin{abstract}
Since its beginnings, every Cycles and Colourings workshop holds one or two open problem sessions; this document contains the problems (together with notes regarding the current state of the art and related bibliography) presented by participants of the \nth{33} edition of the workshop which took place in Nový Smokovec, Slovakia during August \nth{31} -- September \nth{5}, 2025 (see the workshop webpage \url{https://candc.upjs.sk}).
    
\end{abstract}

\tableofcontents

\newpage
\section[The Ramsey question finding a non-nested matching (János Barát)]{Non-nested matching}

\begin{flushleft}
    \textsc{\Large{János Barát}}\smallskip\\
    Alfréd Rényi Institute of Mathematics, Budapest, Hungary\\
    University of Pannonia, Veszprém, Hungary\\
    {\tt barat@renyi.hu}
\end{flushleft}

\bigskip

The following is a classical Ramsey-type question. Let us 2-color the edges of the complete graph and find the largest monochromatic matching.
It is known that we can 2-color the edges of a complete graph on $3k-2$ vertices such that the largest monochromatic matching has $k-1$ edges. If we have one more vertex, then there is always a monochromatic matching of size $k$. 

We change this classical problem as follows.
Let us {\bf order} the vertices of a complete graph $K_m$ and identify the vertices with the set 
$\{1,2,\dots,m\}$. 
Two independent edges $ij$ and $st$ are  {\bf nested}, if $i<s<t<j$ or $s<i<j<t$. If we 2-color the edges, then we would like to find a large monochromatic non-nested matching.

\bigskip
\includegraphics[width=0.8\linewidth]{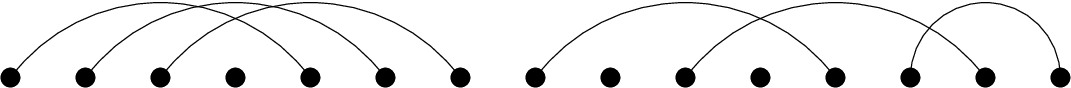}
\bigskip

\begin{boxedproblem}\label{barat1}
    Determine the smallest $m$ such that in every $2$-coloring of the edges of the ordered complete graph $K_m$, there is a monochromatic non-nested matching of size $n$.
\end{boxedproblem}

The answer is trivially between $3n-1$ and $4n-2$. We conjectured that the lower bound is optimal \cite{bgyt}.
A group of young Hungarian researchers (including Benedek Kovács and Kristóf Zólomy) improved the upper bound by a small constant, 5 say. They also checked the small cases up to $n=5$ or 6. 
Recently, we came up with a more substantial improvement of the upper bound. This relies on a Turán-type conjecture, which we can prove with a small error term \cite{bft}. Once this error term is erased, the improvement gives the following.

\begin{boxedproposition}\label{problem2}
   If $m>(2+\sqrt{3})n$, then every $2$-coloring of the edges of an ordered $K_m$ contains a monochromatic non-nested matching with $n$ edges.
\end{boxedproposition}

This gives us more confidence that Problem~\ref{barat1} can be tackled.
We encourage everyone to think about this nice question.

\newpage
\section[Alon-Tarsi number of random 4-regular graphs (Zdeněk Dvořák)]{Alon-Tarsi number of random 4-regular graphs}

\begin{flushleft}
    \textsc{\Large{Zdeněk Dvořák}}\smallskip\\
    Charles University, Prague, Czech Republic\\
    {\tt rakdver@iuuk.mff.cuni.cz}
\end{flushleft}

\bigskip

Let us consider a random 4-regular graph $G$, defined as the union of two uniformly random cycles on the same vertex
set\footnote{Strictly speaking, this is a multigraph; but $G$ is simple with constant probability, and thus this does not affect
the problem.  Moreover, one could think that considering uniformly random 4-regular graphs (say arising from the configuration model) would be more natural.
However, a uniformly random 4-regular graph a.a.s. has two disjoint Hamiltonian cycles~\cite{twoham}, and thus this again does not
play any role in the proposed problem. Another a.a.s. equivalent model one could consider is the union of four uniformly random perfect matchings.}.
Such a graph $G$ is a.a.s. 3-colorable~\cite{threecol}, but it is an open question whether $G$ is
also a.a.s. \emph{3-choosable}, i.e., properly colorable for all assignments of lists of three allowed colors to vertices.

A natural approach to resolving this is by using the Alon-Tarsi method~\cite{alontarsi}:  Consider the \emph{graph
polynomial} $P_G=\prod_{uv\in E(G)} (u-v)$ of a graph $G$, whose variables correspond to the vertices of $G$.
\begin{boxedtheorem}[Alon and Tarsi~\cite{alontarsi}]
Let $G$ be a graph, let $k$ be a positive integer and let $f:V(G)\to \{0,\ldots,k-1\}$ be any function.
If the coefficient of $P_G$ at the monomial $\prod_{v\in V(G)} v^{f(v)}$ is non-zero, then $G$ is $k$-choosable.
\end{boxedtheorem}
Note that the graph polynomial $P_G$ is uniform and every monomial which appears in it with a non-zero coefficient
has the same total degree, equal to $|E(G)|$.  Hence, in our (4-regular) setting, the only relevant monomial is $\prod_{v\in V(G)} v^2$.

\begin{boxedproblem}\label{dvorak1}
    Let $G$ be a random 4-regular graph.  Is it true that the coefficient of $P_G$ at $\prod_{v\in V(G)} v^2$ is a.a.s. non-zero?
\end{boxedproblem}

It is relatively easy to check that (in the two Hamiltonian cycles model) the expected value of the coefficient is zero.
Hence, to prove the claim, one would need some anti-concentration inequality.  Conversely, to disprove it, one would
need to show that the coefficient is tightly concentrated on its expected value.

\newpage
\section[Strong colouring for maximum degree 2 (Penny Haxell)]{Strong colouring for maximum degree 2}

\begin{flushleft}
    \textsc{\Large{Penny Haxell}}\smallskip\\
    University of Waterloo, Waterloo ON, Canada\\
        {\tt pehaxell@uwaterloo.ca}
\end{flushleft}

\bigskip

%Concise introduction: only necessary definitions, maybe a background (who posed the problem originally, connections to other open problems, etc.) 

\begin{boxedproblem}\label{haxell1}
    Let $G$ be a disjoint union of cycles. Let $H$ be a disjoint union of $K_4$'s. Is $G\cup H$ 4-colourable?
\end{boxedproblem}

%\begin{boxedproblem}\label{problem2}
%    Another problem or a question.
%\end{boxedproblem}

If $H$ is a disjoint union of $K_5$'s, then $G\cup H$ is 5-colourable. If $H$ is a disjoint union of $K_3$'s, then $G\cup H$ is not necessarily 3-colourable.

The problem is a special case of a general question about strong chromatic number of graphs: suppose $G$ has maximum degree $d$. If $H$ is a disjoint union of $K_{2d}$'s then is $G\cup H$ ($2d$)-colourable? If this is true, the quantity $2d$ is best possible. It is known that if $H$ is a disjoint union of $K_{3d-1}$'s then $G\cup H$ is $(3d-1)$-colourable.  

See the reference given, and other works that are noted in its list of references or that cite it.
%Short remarks on Problems~\ref {haxell1} and~\ref {problem2}, along with a brief list of relevant references; it can serve as a good starting point for considering how to solve the problem~\cite{ref1,ref2}. 

\newpage
\section[Strong edge coloring of claw-free cubic graphs (František Kardoš)]{Strong edge coloring of claw-free cubic graphs}

\begin{flushleft}
    \textsc{\Large{František Kardoš}}\smallskip\\
    LaBRI, Université de Bordeaux, France\\
    {\tt frantisek.kardos@u-bordeaux.fr}
\end{flushleft}

\bigskip

A strong edge coloring assigns colors to the edges of a given graph in a way that any two edges that are either adjacent to each other or adjacent to a common edge receive distinct colors. In other words, it is a proper (vertex) coloring of the square of the line graph of a given graph. The smallest $k$ such that a given graph $G$ has a strong edge coloring with $k$ colors is called the strong chromatic index of $G$ and is denoted by $\chi'_s(G)$.

Lv, Li, and Zhang \cite{LLZ} proved that if $G$ is a claw-free subcubic graph other than the triangular prism, then $\chi'_s(G)\le 8$. This upper bound was improved to 7 by Lin and Lin \cite{LinLin}, who also found an infinite series of graphs for which the bound is tight. However, all these examples contain diamonds. 

It is a well-known fact that diamond-free claw-free cubic graphs can be obtained from cubic graphs by the operation of \emph{truncation} -- replacing each vertex by a triangle. Let us denote $T(G)$ the cubic graph obtained by truncating a cubic graph $G$.

\begin{boxedproblem}\label{kardos1}
    Is it true that every diamond-free claw-free cubic graph is strongly 6-edge-colorable?
    In other words, is it true that $\chi'_s(T(G))=6$ for every cubic graph $G$?  
\end{boxedproblem}

There is a partial solution by Han and Cui \cite{HanCui} who proved that the truncated prisms are indeed strongly 6-edge-colorable.

\newpage

\section[Homogeneous coloring of cubic graphs (Borut Lu\v{z}ar, Roman Sot\'ak)]{Homogeneous coloring of cubic graphs}

\begin{flushleft}
    \textsc{\Large{Borut Lu\v{z}ar}}\smallskip\\
    Faculty of Information Studies in Novo mesto, Slovenia\\
    Rudolfovo - Science and Technology Centre, Novo mesto, Slovenia.\\
    {\tt borut.luzar@gmail.com}
\end{flushleft}

\begin{flushleft}
    \textsc{\Large{Roman Sot\'ak}}\smallskip\\
    P. J. \v Saf\'arik University in Košice, Slovakia\\
    {\tt roman.sotak@upjs.sk}
\end{flushleft}

\bigskip

A proper coloring of vertices where in the open neighborhood of every vertex exactly $k$ colors appear
is $k$-homogeneous coloring~\cite{JanMad15}.
Clearly, a graph admits $1$-homogeneous coloring if and only if it is bipartite. 
Much more interesting is $2$-homogeneous coloring.
For bipartite cubic graphs, it is known that they admit a $2$-homogeneous coloring with at most $6$ colors~\cite{JanLuzMadaSot17}.
On the other hand, there are cubic graphs that do not admit a homogeneous $2$-coloring;
e.g., consider the graph containing as a subgraph a triangular prism with one edge of a triangle subdivided.

Two questions naturally arise; namely, which cubic graphs admit a $2$-homogeneous coloring,
and what is the minimum number of colors sufficient to color any cubic graph that admits a $2$-homogeneous coloring.
Based on some preliminary research%with R. Sot\'{a}k
, we propose the following two conjectures.

\begin{boxedproblem}\label{luzar1}
	If a cubic graph admits a $2$-homogeneous coloring,
	then $4$ colors are always sufficient.
\end{boxedproblem}

\begin{boxedproblem}\label{luzar2}
    There is a finite number of connected bridgeless cubic graphs that do not admit a $2$-homogeneous coloring. 
\end{boxedproblem}

Note that the first problem appeared in~\cite{JanLuzMadaSot17} in a weaker version, assuming only bipartite cubic graphs.

\newpage
\section[A conjecture on a pair of two flows (Jozef Rajník)]{A conjecture on a pair of two flows}

\begin{flushleft}
    \textsc{\Large{Jozef Rajník}}\smallskip\\
    Comenius University, Bratislava, Slovakia\\
    {\tt jozef.rajnik@fmph.uniba.sk}
\end{flushleft}

\bigskip

\begin{boxedconjecture}
	Each bridgeless graph admits a $2$-flow $\varphi_2$ and a $4$-flow $\varphi_4$ such that for any edge $e$, if $\varphi_2(e) = 0$, then $|\varphi_4(e)| \ge 2$.
\end{boxedconjecture}

This conjecture implies the $5$-flow conjecture since $(5\varphi_2 + \varphi_4)/2$ is a circular $5$-flow. Thus, it offers a possible approach to the $5$-flow conjecture. It can also be viewed as a generalisation of Seymour's approach in the proof of the $6$-flow theorem, where considering a pair of $2$- and $3$-flows played a crucial role. Beside forbidding the pair of flow values $(0, 0)$, our setting also forbids pairs $(0, \pm 1)$. We call this pair of flows a \emph{$1/2$-flow-pair}. We have verified this conjecture for all cyclically $4$-edge-connected snarks up to $34$ vertices, for all snarks with circular flow number $5$ on 36 vertices, and for snarks of oddness $4$ up to $44$ vertices. Further information on such flow pairs will be available in the upcoming manuscript~\cite{Gaborik}.

\newpage
\section[Almost-$\mathbb{Z}_2^2$ and almost-$\mathbb{Z}_4$ connectivities for snarks (Nikolay Ulyanov)]{Almost-$\mathbb{Z}_2^2$ and almost-$\mathbb{Z}_4$ connectivities for snarks}

\begin{flushleft}
    \textsc{\Large{Nikolay Ulyanov}}\smallskip\\
    Berlin, Germany\\
    {\tt ulyanick@gmail.com}
\end{flushleft}

\bigskip

We borrow the terminology from F.~Jaeger et al. \cite{Jaeger}. Let $G=(V,E)$ be a directed graph, and $f$ a mapping from $E$ into a non-trivial Abelian group $A$. Associate with $f$ its \textit{boundary} $\partial f$, a mapping from $V$ to $A$, defined by $\partial f(x) = \sum_{e \text{ leaving } x}f(e) - \sum_{e \text{ entering } x}f(e)$. We say that G is $A$-\textit{connected} if for every $b\colon V \rightarrow A$ with $\sum_{x \in V}b(x)=0$ there is an $f\colon E \rightarrow A - \{0\}$ with $b = \partial f$ (and we call all these boundaries $admissible$).

An $A$-\textit{nowhere-zero-flow} in $G$ is such an $f$ with $\partial f = 0$. W.~T.~Tutte \cite{Tutte} famously conjectured that every bridgeless graph has a $\mathbb{Z}_5$-nowhere-zero-flow. The conjecture could be reduced to considering only cubic graphs, which could be further divided into 2 classes, based on whether they are 3-edge-colourable or not.

Class 2 cubic graphs are not 3-edge-colourable, and it's easy to show that they have neither a $\mathbb{Z}_4$-nowhere-zero-flow nor a $\mathbb{Z}_2^2$-nowhere-zero-flow. In the context of group connectivity, a natural question can be asked about whether other boundaries are admissible or not.

The smallest class 2 cubic graph is the Petersen graph. Interestingly enough, by exhaustive search it's possible to show that for the Petersen graph all boundaries are admissible except for $\partial f = 0$, both for $\mathbb{Z}_2^2$ and $\mathbb{Z}_4$ groups. We will call such graphs \textit{almost}-$A$-\textit{connected} graphs.

Looking at bigger examples of class 2 cubic graphs, a natural subfamily of graphs emerges known as \textit{snarks}, which we define as cubic non-3-edge-colourable cyclically 4-edge-connected graphs with girth at least 5. Looking at the database of snarks, hosted by House of Graphs \cite{Coolsaet}, we'd like to propose the following open problem:

\begin{boxedproblem}
Every bicritical snark is almost-$\mathbb{Z}_2^2$-connected and almost-$\mathbb{Z}_4$-connected.
\end{boxedproblem}

A snark is \textit{bicritical} \cite{NedelaSkoviera} if the removal of any two distinct vertices produces a 3-edge-colourable graph.

As mentioned above, we have verified the proposal for the Petersen graph, and we have also verified it computationally for both Blanuša snarks on 18 vertices. Unfortunately, checking bigger snarks with exhaustive search is quite time-consuming.

The proposed open problem also seems to work "in the opposite direction". We have verified that Tietze's graph, all snarks on 24 vertices (which are all non-bicritical), and all smallest critical non-bicritical snarks (a snark is \textit{critical} if the removal of any two \textit{adjacent} vertices produces a 3-edge-colourable graph) on 32 vertices fail to have "almost-connectivities" mentioned in the proposal.


\begin{thebibliography}{1}

\bibitem{bft} Barát, J., Freschi, A., Tóth, G.
Matchings avoiding ordered patterns.
{\em Manuscript}.

\bibitem{bgyt}
Barát, J., Gyárfás, A., Tóth, G. (2024). Monochromatic spanning trees and matchings in ordered complete graphs. \textit{Journal of Graph Theory}, 105 (4), 523--541.
% J. Barát, A. Gyárfás, G. Tóth.
% Monochromatic spanning trees and matchings in ordered complete graphs.
% Journal of Graph Theory 105 (4), 523--541.
%{\em manuscript}


\end{thebibliography}

\begin{thebibliography}{1}

\bibitem{alontarsi}
Alon, N., Tarsi, M. (1992). Colorings and orientations of graphs. \textit{Combinatorica}, 12 (2), 125--134.
% Noga Alon and Michael Tarsi, \emph{Colorings and orientations of graphs}, Combinatorica 12 (1992), 125--134.

\bibitem{twoham}
Kim, J. H., Wormald, N. C. (2001). Random matchings which induce Hamilton cycles and Hamiltonian decompositions of random regular graphs. \textit{Journal of Combinatorial Theory, Series B}, 81 (1), 20--44.
% Jeong Han Kim and Nicholas C. Wormald, \emph{Random matchings which induce Hamilton cycles and Hamiltonian decompositions of random regular graphs},
% Journal of Combinatorial Theory, Series B 81 (2001), 20--44.

\bibitem{threecol}
Shi, L., Wormald, N. (2007). Colouring random 4-regular graphs. \textit{Combinatorics, Probability and Computing}, 16 (2), 309--344.
% Lingsheng Shi and Nicholas C. Wormald, \emph{Colouring random 4-regular graphs}, Combinatorics, Probability and Computing 16 (2007), 309--344.

\end{thebibliography}

\begin{thebibliography}{1}

\bibitem{ref1}
Haxell, P. E. (2004). On the strong chromatic number. \textit{Combinatorics, Probability and Computing}, 13 (6), 857--865.
% Haxell, P.E., On the strong chromatic number, {\it Combinatorics, Probability and Computing} 13 (2004),  857--865

\end{thebibliography}

\begin{thebibliography}{1}

\bibitem{HanCui}
Han, Z., Cui, Q. (2023). A note on strong edge-coloring of claw-free cubic graphs. \textit{Journal of Applied Mathematics and Computing}, 69 (3), 2503--2508.
% Zhenmeng Han and Qing Cui, A note on strong edge-coloring of claw-free cubic graphs, Journal of Applied Mathematics and Computing 69 (2023), 2503--2508.
% https://link.springer.com/article/10.1007/s12190-023-01847-x

\bibitem{LinLin}
Lin, Y., Lin, W. (2023). The tight bound for the strong chromatic indices of claw-free subcubic graphs. \textit{Graphs and Combinatorics}, 39 (3), 58.
% Yuquan Lin and Wensong Lin, The Tight Bound for the Strong Chromatic Indices of Claw-Free Subcubic Graphs Graphs and Combinatorics 39, 58 (2023).
% https://link.springer.com/article/10.1007/s00373-023-02655-7

\bibitem{LLZ}
Lv, J. B., Li, J., Zhang, X. (2022). On strong edge-coloring of claw-free subcubic graphs. \textit{Graphs and Combinatorics}, 38 (3), 63.
% Jian-Bo Lv, Jianxi Li, and Xiaoxia Zhang, On Strong Edge-Coloring of Claw-Free Subcubic Graphs, Graphs and Combinatorics 38, 63 (2022). 
% https://link.springer.com/article/10.1007/s00373-022-02462-6


\end{thebibliography}

\begin{thebibliography}{1}

\bibitem{JanMad15}
Janicov\'{a}, M. (2015). Okolia vrcholov v zafarben\'{y}ch grafoch. \textit{Master's thesis}, P.~J.~\v{S}af\'arik University in Ko\v{s}ice.

\bibitem{JanLuzMadaSot17}
Janicová, M., Madaras, T., Soták, R., Lužar, B. (2017). From NMNR-coloring of hypergraphs to homogenous coloring of graphs. \textit{Ars Mathematica Contemporanea}, 12 (2), 351--360.
% M. Janicov\'{a}, B. Lu\v{z}ar, T. Madaras, R. Sot\'{a}k. (2017). From NMNR-coloring of hypergraphs to homogeneous coloring of graphs. \textit{}, \textit{12}, 351--360.


\end{thebibliography}

\begin{thebibliography}{1}

\bibitem{Gaborik}
Gáborik, L., Kurz, S., Mazzuoccolo, G., Rajník, J., Rieg, F. Manhattan and Chebyshev flows. \textit{Manuscript}.
% L. Gáborik, S. Kurz, G. Mazzuoccolo, J. Rajník, F. Rieg, \emph{Manhattan and Chebyshev flows}, Manuscript.

\bibitem{Seymour} 
Seymour, P. D. (1981). Nowhere-zero 6-flows. J\textit{ournal of combinatorial theory, series B}, 30 (2), 130--135.
% P. D. Seymour, \emph{Nowhere-zero 6-flows}, J. Comb. Theory Ser. B, 30(2), pp. 130--135, 1981.

\end{thebibliography}

\begin{thebibliography}{1}

\bibitem{Coolsaet}
Coolsaet, K., D’hondt, S., Goedgebeur, J. (2023). House of Graphs 2.0: A database of interesting graphs and more. \textit{Discrete Applied Mathematics}, 325, 97--107. Available at \url{https://houseofgraphs.org}.
% Coolsaet, K., D’hondt, S., and Goedgebeur, J.  
% \newblock House of Graphs 2.0: A database of interesting graphs and more.  
% \newblock {\em Discrete Applied Mathematics}, 325 (2023), 97--107.  
% \newblock Available at \url{https://houseofgraphs.org/meta-directory/snarks}.

\bibitem{Jaeger}
Jaeger, F., Linial, N., Payan, C., Tarsi, M. (1992). Group connectivity of graphs—a nonhomogeneous analogue of nowhere-zero flow properties. \textit{Journal of Combinatorial Theory, Series B}, 56 (2), 165--182.
% Jaeger, F., Linial, N., Payan, C., and Tarsi, M.  
% \newblock Group connectivity of graphs—a nonhomogeneous analogue of nowhere-zero flow properties.  
% \newblock {\em Journal of Combinatorial Theory, Series B}, 56(2) (1992), 165--182.

\bibitem{NedelaSkoviera}
Nedela, R., \v Skoviera, M. (1996). Decompositions and reductions of snarks. \textit{Journal of Graph Theory}, 22 (3), 253--279.
% Nedela, R., and Škoviera, M.  
% \newblock Decompositions and reductions of snarks.  
% \newblock {\em Journal of Graph Theory}, 22(3) (1996), 253--279.

\bibitem{Tutte}
Tutte, W. T. (1954). A contribution to the theory of chromatic polynomials. \textit{Canadian journal of mathematics}, 6, 80--91.
% Tutte, W. T.  
% \newblock A contribution to the theory of chromatic polynomials.  
% \newblock {\em Canadian Journal of Mathematics}, 6 (1954), 80--91.





\end{thebibliography}
\end{document}